\documentclass[preprint,12pt]{elsarticle}
\usepackage{amssymb}
\usepackage{amsmath}
\usepackage{amssymb}
\usepackage{enumerate}
\usepackage{hyperref}
\usepackage{amsthm}
\usepackage[utf8]{inputenc}
\usepackage[T1]{fontenc}
\newtheorem{theorem}{Theorem}
\newtheorem{definition}[theorem]{Definition} 
\newtheorem{lemma}{Lemma} 
\newtheorem{remark}[theorem]{Remark}
\newtheorem{proposition}[theorem]{Proposition}

\journal{nothing}

\begin{document}
\begin{frontmatter}
\title{Regularity of solutions to time-harmonic Maxwell’s system with
	Hölder and various lower than Hölder continuous coefficients}

\author[add1]{Lei Yu } %% Author name
\author[add2]{Basang Tsering-xiao\corref{cor1}}

\address[add1]{Mathematical Department of Tibet University, Lhasa, Tibet, 850000, PR China} 
\address[add2]{Mathematical Department of Tibet University, Lhasa, Tibet, 850000, PR China} 
\cortext[cor1]{Corresponding author\\
E-mail addresses:  basangtu@qq.com (B. Tsering-xiao), leiyu@stu.utibet.edu.cn(L.Yu)}
\begin{abstract}
     The purpose of this paper is to establish a complete Schauder theory for the second-order linear elliptic equation and the time-harmonic Maxwell's system. We prove global Hölder regularity for the solutions to the time-harmonic anisotropic Maxwell's equations under Hölder continuous coefficients, raising the Hölder index to the interval \( (0, 1) \). We also show that when the coefficients are \( C^0 \), we can establish \( L^p \) theory for the solutions. When one of the coefficients is \( L^\infty \), we obtain Hölder regularity results for the electric field.
     
\end{abstract}
\end{frontmatter}
\section{Introduction}
Weak solutions  $E, H \in H(\operatorname{curl},\Omega ) := \{F \in
L^2(\Omega ) : \operatorname{curl} F \in L^2(\Omega )\}$,
\begin{equation}\label{1}
\begin{cases}\operatorname{curl}H=\text{i}\omega\varepsilon E+J_e&\text{ in }\Omega,\\\operatorname{curl}E=-\text{i}\omega\mu H+J_m&\text{ in }\Omega,\\E\times\nu=G\times\nu&\text{ on }\partial\Omega,\end{cases}
\end{equation}
where $\Omega\subseteq\mathbb{R}^3$ be an open bounded domain with $C^
{1,1}$ boundary,this boundary condition holds when $\Omega$ is
surrounded by a perfect conductor, see for instance \cite{b1,b2,peter} and references therein. the terms $J_m,J_e \in L^2(\Omega)$
 are current sources
and the boundary condition G is in $H(curl,\Omega)$, where $\varepsilon$ and $\mu $ are real matrix-valued functions in $L^\infty \left(\mathbb{R}^3\right)$,
such that for all $\xi\in\mathbb{R}^3$ and almost every $x\in\mathbb{R}^3$, there exist $0 < \lambda _1 \leq \lambda _2 < \infty$ and
$0 < \lambda _3 \leq \lambda _4 < \infty$.
\begin{equation}\label{2}
\lambda_1|\xi|^2\leq\varepsilon(x)\xi\cdot\xi\leq\lambda_2|\xi|^2, \\\lambda_3|\xi|^2\leq\mu(x)\xi\cdot\xi\leq\lambda_4|\xi|^2.
\end{equation}

The aim of this paper is to investigate the regularity of electromagnetic fields under low-regularity assumptions concerning material parameters. The material parameters \((\varepsilon(x), \mu(x))\) represent the matrix-valued coefficients in Maxwell's equations.

This research is primarily driven by the investigation of electromagnetic inverse problems in liquid crystals, as previously presented by Basang et al. in \cite{Basang2}. Beyond the regularity of coefficients and outcomes, a plethora of issues related to Maxwell's equations has been explored, including boundary regularity and the extension of uniqueness \cite{W,W2,W3,W4}. Several scholars, such as those cited in \cite{C, J, JM, Yin, Peter}, have examined the regularity of solutions to Maxwell's equations, particularly in scenarios characterized by high-regularity coefficients.

Up to the present, Alberti, Capdeboscq \cite{A3,A2,A1}, and Basang \cite{Basang1} represent a very limited number of researchers who have addressed the regularity of coefficients that are less than Lipschitz continuous. We assume Hölder continuity for \( \mu \) and weaken the regularity of \( \varepsilon \). Due to the symmetry of the equation, we can derive the weakened regularity of \( \varepsilon \) when \( \mu \) is Hölder continuous, which leads to the regularity of \( H \). However, since the boundary conditions pertain to \( E \), applying symmetry results in different situations.

This paper extends the regularity results presented in \cite{A3}, which demonstrated that one can obtain \( W^{1,6} \) from \( L^2 \) using classical elliptic regularity and Sobolev embedding, leading to the result of \( C^{0,1/2} \). We find that the interval \( (0, 1/2] \) is not perfect.

By adjusting the parameter \( p \) in \( L^p \), it is possible to achieve any \( C^{0,\alpha} \). However, this technique requires refinement and necessitates additional iterations, while meticulously controlling the parameter \( p \).

This endeavor represents a significant step toward constructing a comprehensive regularity theory, which is a crucial advancement in the field. Below, we provide a historical review of the regularity of Maxwell's equations.

\begin{itemize}
	\item 1955, Friedrichs: if $u\in L^2,\nabla \times u \in L^2,\mathrm{div} u\in L^2$ then $u\in W^{1,2}$. 
	\item 1981, Weber: if $\varepsilon,\mu \in W^{1,\infty}$ then $E, H \in W^{1,2}$. 
	\item 2004, Yin: if $\varepsilon \in W^{1,\infty} $ then $E \in C^{0,\alpha}$. 
	\item 2014, Yves Capdeboscq: if $\varepsilon,\mu \in W^{1,3+\delta}$, then $E, H \in W^{1,2}$. 
	\item2018, Alberti : $\varepsilon,\mu \in C^{0,\alpha}$, then $E, H \in C^{0,\alpha}$  ($\alpha \in (0,1/2]$). 
	\item2019, Basang: $\varepsilon,\mu \in W^{1,3+\delta}$, and $\|\mu(x)-I\|_{\mathcal{L}^{\infty}(\Omega)^{3\times3}}+\|\varepsilon(x)-I\|_{\mathcal{L}^{\infty}(\Omega)^{3\times3}}<\delta_{0}$, then $E, H \in C^{1,\alpha}$  . 
	\item2024 our work $\varepsilon,\mu \in C^{0,\alpha}$, then $E, H \in C^{0,\alpha}$  ($\alpha \in (0,1)$). 
\\We imporve Hölder index $1/2$ to $1$ from 2018 Alberti 
\item 2024 our work $\varepsilon \in C^{0,\alpha}$, $\mu \in C^0, $ then $E \in C^{0,\alpha},H \in \mathcal{L}^{2,3-\tau}$  ($\alpha \in (0,1)$) and for any $\tau > 0$. 
\item 2024 our work $\varepsilon,\mu \in C^{0}$, then $E,H \in \mathcal{L}^{2,3-\tau}$ and for any $\tau > 0$. 
\item 2024 our work $\varepsilon \in C^{0,\alpha},\mu \in L^\infty$, then $E \in C^{0,\alpha}$ ($\alpha \in (0,\delta_{0}]$).

Of course, the work of the harmonic Maxwell equations is not limited to this, and we will only list the work where the regularity of the coefficients is very low

\end{itemize}
\section{ Main results}
\begin{theorem}\label{m}
Assume that \eqref{2} holds true and  $\varepsilon(x),\mu(x) \in C^{0,\alpha} (\bar{\Omega},\mathbb{R}^3)$,
for some $\alpha \in (0,1)$, take $J_{e}, J_{m}\in C^{0,\alpha}(\overline{\Omega},\mathbb{R}^3 $) and $G\in C^{1,\alpha}(\operatorname{curl},\Omega)$, Let $(E,H)\in H(curl,\Omega)$  be a weak solution of \eqref{1}.  Then $E, H\in C^{0,\alpha}(\overline{\Omega})$ and
$$\|(E,H)\|_{C^{0,\alpha}(\overline{\Omega})}\leq C\big(\|(E,H)\|_{L^2(\Omega)}+\|G\|_{C^{1,\alpha}(curl,\Omega)}+\|(J_e,J_m)\|_{C^{0,\alpha}(\overline{\Omega})}\big),$$
for some constant C depending only on $\Omega$, $\left\|\varepsilon\right\|_{C^{0,\alpha}\left(\overline{\Omega}\right)}$, $\left\|\mu\right\|_{C^{0,\alpha}\left(\overline{\Omega}\right)}$ and  $\omega$. 
\end{theorem}
\begin{theorem}\label{b}
	Assume that \eqref{2} holds true and  $\varepsilon(x),\mu(x) \in C^{0} (\bar{\Omega},\mathbb{R}^3)$, take $J_{e}, J_{m}\in L^{q}(\overline{\Omega},\mathbb{R}^3 $) and $G\in W^{1,q}(\operatorname{curl},\Omega)$, for  $q>3$, Let $(E,H)\in H(curl,\Omega)$  be a weak solution of \eqref{1} Then $E, H\in L^{q}({\Omega})$ and
	$$\|(E,H)\|_{L^{q}(\overline{\Omega})}\leq C\big(\|(E,H)\|_{L^2(\Omega)}+\|G\|_{W^{1,q}(curl,\Omega)}+\|(J_e,J_m)\|_{{L^q}(\overline{\Omega})}\big),$$
	for some constant C depending only on $\Omega$, $\left\|\varepsilon\right\|_{C^{0}\left(\overline{\Omega}\right)}$ $\left\|\mu\right\|_{C^{0}\left(\overline{\Omega}\right)}$ and $\omega$.  
\end{theorem}

\begin{remark}
	In Theorem \ref{b}, take $J_{e}, J_{m}\in L^{\infty}(\overline{\Omega},\mathbb{R}^3 $) and  $G\in W^{1,\infty}(\operatorname{curl},\Omega)$, we can not obtain $E, H\in \mathcal{L}^{2,3}({\Omega})$, only obtain $E,H \in \mathcal{L}^{2,3-\tau}({\Omega})$ and for any $\tau > 0$. 
\end{remark}

\begin{theorem}\label{c}
	Assume that \eqref{2} holds true and  $\varepsilon(x), \in C^{0,\alpha} (\bar{\Omega},\mathbb{R}^3),\mu \in C^{0} (\bar{\Omega},\mathbb{R}^3)$,
	for some $\alpha \in (0,1)$, take $G,J_{e}\in C^{0,\alpha}(\overline{\Omega},\mathbb{R}^3 $), $\operatorname{curl}G,J_m\in L^{\infty}(\bar{\Omega},\mathbb{R}^3)$, Let $(E,H)\in H(curl,\Omega)$  be a weak solution of \eqref{1} Then $E \in C^{0,\alpha} (\bar{\Omega})$ and
	$$\|E\|_{C^{0,\alpha}(\overline{\Omega})}\leq C\big(\|(E,H)\|_{L^2(\Omega)}+\|(G,J_e)\|_{C^{0,\alpha}(\overline{\Omega})}+\|\mathrm{(curl}G,J_m)\|_{L^\infty(\overline{\Omega})}\big),$$
	for some constant C depending only on $\Omega$, $\left\|\varepsilon\right\|_{C^{0,\alpha}\left(\overline{\Omega}\right)}$, $\left\|\mu\right\|_{C^{0}\left(\overline{\Omega}\right)}$ and $\omega$.  
	
\end{theorem}

\begin{theorem}\label{d}
	Assume that \eqref{2} holds true and  $\varepsilon(x) \in C^{0,\alpha} (\bar{\Omega},\mathbb{R}^3),\mu \in L^\infty (\bar{\Omega},\mathbb{R}^3)$,
	for some $\alpha \in (0,\delta_{0}]$, take $J_{e},G\in C^{0,\alpha}(\overline{\Omega},\mathbb{R}^3 $), $\operatorname{curl}G,J_m\in \mathcal{L}^{2,\gamma_0}(\bar{\Omega},\mathbb{R}^3)$, Let $(E,H)\in H(curl,\Omega)$  be a weak solution of \eqref{1} Then $E \in C^{0,\alpha} (\bar{\Omega})$ and
	$$\|E\|_{C^{0,\alpha}(\overline{\Omega})}\leq C\big(\|(E,H)\|_{L^2(\Omega)}+\|(G,J_e)\|_{C^{0,\alpha}(\overline{\Omega})}+\|(\operatorname{curl}G,J_m)\|_{\mathcal{L}^{2,\gamma_0}(\overline{\Omega})}\big),$$
	for some constant C depending only on $\Omega$, $\left\|\varepsilon\right\|_{C^{0,\alpha}\left(\overline{\Omega}\right)}$,$\left\|\mu\right\|_{C^{0}\left(\overline{\Omega}\right)}$ and $\omega$. \\ $\delta_0,\gamma_0$ define from Lemma 3,4
	
\end{theorem}

\begin{remark}
	In  Theorem \ref{d}.  
	We improve the results in \cite{A3} such that $\delta_0$  is not affected by 1/2.

\end{remark}

Let us consider the homogeneous Maxwell's equations with homogeneous boundary condition . 
\begin{equation}\label{3}
	\begin{cases}\operatorname{curl}H=\text{i}\omega\varepsilon E&\text{ in }\Omega,\\\operatorname{curl}E=-\text{i}\omega\mu H&\text{ in }\Omega,\\E\times\nu=0&\text{ on }\partial\Omega,\end{cases}
\end{equation}

\begin{proposition}\label{4}
Assume that \eqref{2} holds true and  $\varepsilon(x),\mu(x) \in C^{0,\alpha} (\bar{\Omega},\mathbb{R}^3)$,
for some $\alpha \in (0,{1}) $.  Let $(E,H)\in H(curl,\Omega)$  be a weak solution of \eqref{3} Then $E, H\in C^{0,\alpha}(\overline{\Omega})$ and
$$\|(E,H)\|_{C^{0,\alpha}(\overline{\Omega})}\leq C\|(E,H)\|_{L^2(\Omega)},$$
for some constant C depending only on $\Omega$, $\left\|\varepsilon\right\|_{C^{0,\alpha}\left({\Omega}\right)}$,$\left\|\mu\right\|_{C^{0,\alpha}\left(\overline{\Omega}\right)}$ and $\omega$. 
\end{proposition}

\begin{proposition}\label{l}
	Assume that \eqref{2} holds true and  $\varepsilon(x) \in C^{0,\alpha}\left(\mathbb{R}^3\right),\mu \in  C^0\left(\mathbb{R}^3\right) $,
	for some $\alpha \in(0,1)$,
	Let $(E,H)\in H(curl,\Omega)$  be a weak solution of \eqref{3} Then $E\in C^{0,\alpha}(\overline{\Omega})$ and
	$$\|E\|_{C^{0,\alpha}(\overline{\Omega})}\leq C\|E\|_{L^2(\Omega)},$$
	for some constant C and   depending only on $\Omega$, $\left\|\varepsilon\right\|_{C^{0,\alpha}\left(\overline{\Omega}\right)}$,
	$\left\|\mu \right\|_{C^0(\overline{\Omega})}$ and $\omega$. \\
\end{proposition}

\begin{proposition}\label{f}
	Assume that \eqref{2} holds true and  $\varepsilon(x) \in C^{0,\alpha}\left(\mathbb{R}^3\right),\mu \in  L^\infty \left(\mathbb{R}^3\right)$,
	for some $\alpha \in(0,\delta_0]$,
	Let $(E,H)\in H(curl,\Omega)$  be a weak solution of \eqref{3} Then $E\in C^{0,\alpha}(\overline{\Omega})$ and
	$$\|E\|_{C^{0,\alpha}(\overline{\Omega})}\leq C\|E\|_{L^2(\Omega)},$$
	for some constant C and $\delta_0$ (both independent of $(E,H)$) depending only on $\Omega$, $\left\|\varepsilon\right\|_{C^{0,\alpha}\left(\overline{\Omega}\right)}$,
	$\left\|\mu \right\|_{L^\infty\left(\overline{\Omega}\right)}$ and $\omega$.  $\delta_0$ define from Lemma 3.\\
\end{proposition}

\begin{proposition}\label{g}
	Assume that \eqref{2} holds true and  $\varepsilon(x),\mu(x) \in C^{0}\left(\mathbb{R}^3\right)$,Let $(E,H)\in H(curl,\Omega)$  be a weak solution of \eqref{3} Then $E,H \in \mathcal{L}^{2,3-\tau}({\Omega})$ and for any $\tau > 0$ .

	$$\|(E,H)\|_{\mathcal{L}^{2,3-\tau}(\overline{\Omega})}\leq C\|(E,H)\|_{L^2(\Omega)},$$
	for some constant C  depending only on $\Omega$, $\left\|\varepsilon\right\|_{C^{0}\left(\overline{\Omega}\right)}$,
	$\left\|\mu \right\|_{C^{0}\left(\overline{\Omega}\right)}$. 
\end{proposition}

In summary, proposition \ref{4},\ref{l},\ref{f},\ref{g} is the homogeneous boundary version of Theorem \ref{m},\ref{b},\ref{c},\ref{d} for the homogeneous equation.

\section{Preliminary results}

 In this paper, we are interested in the low regularity of the coefficients.  We do not consider the regularity issues of the domain and boundary, we always assume that the domain and boundary are simply connected.

Assume
$$\lambda|\xi|^2\leq A(x)\xi\cdot\xi\leq\Lambda|\xi|^2. $$ 
with  constants $0 < \lambda < \Lambda< \infty$.

\begin{lemma}[$L^p$ estimate (\cite{p}, Theorem 1)]
	Assume $A(x)\in C^0({\Omega}),\Omega \in C^{1,1} $, then $\operatorname{div} (A(x) \nabla u) =div f$, for $p\in (1,\infty)$, $$f\in L^{P}({\Omega}) \Rightarrow \nabla u\in L^{P}({\Omega}). $$
\end{lemma}

\begin{lemma}[Schauder estimate(\cite{elliptic2}]
	Assume $A(x)\in C^{k,\alpha},\Omega \in C^{k,\alpha}({\Omega})$, $\operatorname{div} (A(x) \nabla u) =div f$,
	$$f\in C^{k,\alpha}({\Omega}) \Rightarrow \nabla u\in C^{k,\alpha}({\Omega}). $$
\end{lemma}

\begin{lemma}[De Giorgi-Nash theorem(\cite{elliptic}, theorem 2. 14)]

	Assume $A(x)\in L^\infty({\Omega}) $, then $\operatorname{div} (A(x) \nabla u) =0$, then  $$ u \in C^{0,\delta_0}({\Omega}). $$
	
	constant $\delta_0 \in (0, 1)$  depends on  $\| A\|_{L^{\infty}}$, $\frac{\lambda}{\Lambda} $ and $\Omega$. 
\end{lemma}

\begin{lemma}[Campanato estimate(\cite{elliptic}, theorem 2. 19)]
	Assume $A(x)\in L^\infty({\Omega}) $, $\operatorname{div} (A(x) \nabla u) =div f$,  $f \in L^{2,\gamma}$, $0 \leq \gamma \leq \gamma_0$ then $$\nabla u \in L^{2,\gamma}({\Omega}). $$ 
	
	$\gamma_0=N-2+2\delta_0$, $\delta_0$ define from Lemma 3,
	
	in particular, $u \in C^{0,\delta}({\Omega})$  with  $\delta =\frac{\gamma-N+2}{2}  $ if $ \gamma > N - 2$. 
\end{lemma}

\begin{lemma}[ Gaffney-Friedrichs Inequality  (\cite{A1}, lemma 2)]
	Take $p \in (1,\infty)$ and $F \in L^p(\Omega)$, such that $\operatorname{curl} F \in L^p(\Omega)$,
	$\mathrm{div} F \in L^p(\Omega)$ and either $F \cdot \nu = 0$ or $F \times \nu = 0$ on $\partial \Omega$.  Then $F \in W^{1,p}(\Omega) $. 
	 $$\|F\|_{W^{1,p}(\Omega)}\leq C(\left\|\mathrm{div} F\right\|_{L^p(\Omega)}+\left\|\operatorname{curl} F\right\|_{L^p(\Omega)}),$$
	for some $C > 0$ depending only on $\Omega$ and $p$. 
\end{lemma}

\section{ Proof of  results}

\subsection{Proof of Proposition \ref{4}. }
\begin{proof}
	\begin{enumerate}[Step 1:]	
		\item  Helmholtz decomposition.  There exist $q_E\in H_0^1(\Omega),q_H\in H^1(\Omega)\mathrm{~and~}\Phi_E,\Phi_H\in H^1(\Omega)\mathrm{~such~that}\\E=\nabla q_E+\operatorname{curl}\Phi_E\quad\operatorname{in} \Omega,\quad H=\nabla q_H+\operatorname{curl}\Phi_H\quad\operatorname{in} \Omega,$\\

		$\begin{cases}\operatorname{div}\Phi_E=0\text{ in }\Omega,\\\Phi_E\cdot\nu=0\text{ on }\partial\Omega,&\end{cases}\quad\begin{cases}\operatorname{div}\Phi_H=0\text{ in }\Omega,\\\Phi_H\times\nu=0\text{ on }\partial\Omega. &\end{cases}$
		
		$$\|\Phi_E\|_{H^1(\Omega)}\leq C\left\|E\right\|_{L^2(\Omega)},\quad\|\Phi_H\|_{H^1(\Omega)}\leq C\left\|H\right\|_{L^2(\Omega)}. $$
		
		\item $L^{p^*}$ estimate for   $H,E$.

		Set $\Psi_E=\operatorname{curl}\Phi_E$.  Observe that for every test function $\phi \in C^{\infty}$,
		
		since integration by parts  \cite{p}, boundary values in the Sobolev space H$^{-1/2}$
	 defined in the distributional sense
		by the natural extension of the Green formulas. 
		
	$$	\begin{array}{rcl}(\operatorname{curl}\vec{u}, \vec{v}) - (\vec{u},\operatorname{curl}\vec{v})&=&<\nu\times\vec{u}, \vec{v}>,\\(\operatorname{div}\vec{u}, \varphi) + (\vec{u},\operatorname{grad}\varphi)&=&<\nu\cdot\vec{u}, \varphi>,\end{array}$$

		$$
			\int_{\Omega}\operatorname{curl}(\nabla q_E)\cdot\Phi-\nabla q_E\cdot\operatorname{curl}\Phi dx =$$
		$$\int_{\Omega}q_{E}\operatorname{div}(\operatorname{curl}\Phi) dx-\int_{\partial\Omega}q_{E}\operatorname{curl}\Phi\cdot\nu ds \\
			=0, 
		$$
		
		which \cite{peter} implies that $\nabla q_E\times\nu=0\mathrm{~on~}\partial\Omega $, so we obtain 
		$$\Psi_E\times\nu=(\text{curl }\Phi_E)\times\nu=E\times\nu-\nabla q_E\times\nu=0\quad\text{on }\partial\Omega. $$ 
		
		Then $\Psi_E$ satisfy
		\begin{equation}
			\begin{cases}\operatorname{curl}\Psi_E=-\text{i}\omega\mu H&\text{ in }\Omega,\\\mathrm{div}\Psi_E=0&\text{ in }\Omega,\\\Psi_E\times\nu=0&\text{ on }\partial\Omega,\end{cases}
		\end{equation}
		
		We can obtain the result through Lemma 5, divergence, curl, and boundary conditions. 
		
		\begin{equation}\label{5}
			\|\Psi_E\|_{W^{1,p}(\Omega)}\leq C\|H\|_{L^p(\Omega)},
		\end{equation}
		$$\operatorname{curl}\Phi_E= \Psi_E\in {W^{1,p}(\Omega)}. $$
		
		The same method applies to the $\operatorname{curl}\Phi_H$. 
		
		Set $\Psi_H=\operatorname{curl}\Phi_H$.  Observe that for every test function $\phi \in C^{\infty}$,
		$$\int_\Omega\Psi_H\cdot\nabla\varphi+\varphi\mathrm{~div~}\Psi_Hdx=\int_\Omega\mathrm{curl~}\Phi_H\cdot\nabla\varphi dx$$
		$$=\int_\Omega\Phi_H\cdot\mathrm{curl~}\nabla\varphi dx+\int_{\partial\Omega}(\Psi_H \times   \nu) \cdot \nabla \phi
		ds=0,$$
		then   $\Psi_H \cdot \nu =0$ on $\partial\Omega$,
		
		\begin{equation}
			\begin{cases}\operatorname{curl}\Psi_H=\text{i} \omega\varepsilon E&\text{ in }\Omega,\\\mathrm{div}\Psi_E=0&\text{ in }\Omega,\\\Psi_H\cdot\nu=0&\text{ on }\partial\Omega. \end{cases}
		\end{equation}
		
		We have
		\begin{equation}\label{7}
			\|\Psi_H\|_{W^{1,p}(\Omega)}\leq C\|E\|_{L^p(\Omega)},
		\end{equation}
		$$\operatorname{curl}\Phi_H= \Psi_E\in {W^{1,p}(\Omega)}. $$
		Thus, we have obtained two results. 
		
	    $$\operatorname{curl}\Phi_H= \Psi_E\in {W^{1,p}(\Omega)},$$
		$$\operatorname{curl}\Phi_E= \Psi_E\in {W^{1,p}(\Omega)}. $$
		By adjusting the size of \( p \), we control it to the  $(1,\frac{3}{2})$,
		
		$$E,H \in L^2(\Omega) \Rightarrow E,H \in L^p(\Omega),p\in (1,\frac{3}{2}),$$
		then through Sobolev embedding $\operatorname{curl}\Phi_H,\operatorname{curl}\Phi_E \in L^{p^*}(\Omega)$, and $p^{*}=\frac{3p}{3-p}\in(\frac{3}{2},3)<3. $

		We will next use the other two sets of equations to estimate $\nabla q_E$ and $\nabla q_H$. 
		\begin{equation}\label{8}
			\begin{cases}-\mathrm{div}(\varepsilon \nabla q_E)=\mathrm{div}(\varepsilon \operatorname{curl}\Phi_E)&\text{ in }\Omega,\\ \nabla q_E=0&\text{ on }\partial\Omega,\end{cases}
		\end{equation}
		\begin{equation}\label{9}
			\begin{cases}-\mathrm{div}(\mu \nabla q_H)=\mathrm{div}(\mu \operatorname{curl}\Phi_H)&\text{ in }\Omega,\\ -(\mu \nabla q_H)\cdot \nu=(\mu \operatorname{curl}\Phi_H)\cdot \nu&\text{ on }\partial\Omega. \end{cases}
		\end{equation}
		Therefore, by the $L^p$ theory for elliptics (Calderón–Zygmund estimate), we obtain $\nabla q_E,\nabla q_H \in L^{p^*}(\Omega)$,\\so $E=\nabla q_E+\operatorname{curl} \Phi_E \in L^{p^*}(\Omega)$ and $H=\nabla q_H+\operatorname{curl} \Phi_H \in L^{p^*}(\Omega)$
		and
		$$\|(E,H)\|_{L^{p^*}(\Omega)}\leq C\|E,H\|_{L^p(\Omega)}. 
		$$
		
		\item $L^{p^{**}}$ estimate for $E,H$\\
		We repeatedly perform Step 2 and Step3. 
		By \eqref{5} \eqref{7},we obtain  
		$$\operatorname{curl} \Phi_E,\operatorname{curl} \Phi_H \in W^{1,p^{*}}(\Omega) \hookrightarrow L^{p^{**}}(\Omega),$$
		and 
		$$\|(\operatorname{curl} \Phi_E,\operatorname{curl} \Phi_H)\|_{W^{1,p^{*}}(\Omega)}\leq C\|(E,H)\|_{L^p{^{*}}(\Omega)},$$
		by the $L^p$ theory for \eqref{8} \eqref{9} (Calderón–Zygmund estimate), we obtain $\nabla q_E,\nabla q_H \in L^{p^{**}}(\Omega)$. 
		\\so $E=\nabla q_E+\operatorname{curl} \Phi_E \in L^{p^{**}}(\Omega)$ and $H=\nabla q_H+\operatorname{curl} \Phi_H \in L^{p^{**}}(\Omega)$.

		Use sobolev embdding once again.   $p^{**}>3$,We have increased the range of the H$\ddot{o}$lder index, $$\alpha=1-\frac{3}{p^{**}} \in(0,1),$$
		$$\|(\operatorname{curl} \Phi_E,\operatorname{curl} \Phi_H)\|_ {C^{0,\alpha}(\bar{\Omega})}\leq C\|(\operatorname{curl} \Phi_E,\operatorname{curl} \Phi_H)\|_{W^{1,p^{**}}(\Omega)}$$
		$$\leq C\|(E,H)\|_{L^{p^{**}}(\Omega)}. $$
		
		\item$ C^{0,\alpha}$ estimate.  \\
		By applying Schauder esitmate for \eqref{8} \eqref{9},
		$$\nabla q_E,\nabla q_H \in C^{0,\alpha}(\overline{\Omega}),$$
		and
		$$\|\nabla q_E,\nabla q_H \|_{C^{0,\alpha}(\overline{\Omega})} \leq C\|(\operatorname{curl} \Phi_E,\operatorname{curl} \Phi_H) \|_{C^{0,\alpha}(\overline{\Omega})} . $$ \\
		Finally 
		$E=\nabla q_E+\operatorname{curl} \Phi_E \in C^{0,\alpha}(\overline{\Omega})$ and $H=\nabla q_H+\operatorname{curl} \Phi_H \in C^{0,\alpha}(\overline{\Omega}). $
	\end{enumerate}	
	
	\end{proof}

\newpage

\subsection{Proof of Theorem \ref{m}. }
\begin{proof}
	In the non-homogeneous process, the influences of \( G \), \( J_e \), and \( J_m \) will arise. 
	
	\begin{enumerate}[Step 1:]	
		\item  Helmholtz decomposition.  There exist $q_E\in H_0^1(\Omega),q_H\in H^1(\Omega)\mathrm{~and~} \\ \Phi_E,\Phi_H\in H^1(\Omega)\mathrm{~such~that}\\E-G=\nabla q_E+\operatorname{curl}\Phi_E\quad\operatorname{in} \Omega,\quad H=\nabla q_H+\operatorname{curl}\Phi_H\quad\operatorname{in},$\\

		$\begin{cases}\operatorname{div}\Phi_E=0\text{ in }\Omega,\\\Phi_E\cdot\nu=0\text{ on }\partial\Omega,&\end{cases}\quad\begin{cases}\operatorname{div}\Phi_H=0\text{ in }\Omega,\\\Phi_H\times\nu=0\text{ on }\partial\Omega. &\end{cases}$
		
		$$\|\Phi_E\|_{H^1(\Omega)}\leq C\left\|E\right\|_{L^2(\Omega)},\quad\|\Phi_H\|_{H^1(\Omega)}\leq C\left\|H\right\|_{L^2(\Omega)}. $$
		
		\item $L^{p^{*}}$ estimate for E, H .

		Set $\Psi_E=\operatorname{curl}\Phi_E$.  Observe that for every test function $\phi \in C^{\infty}$,
		
		since integration by parts .  \\
		$$
			\int_{\Omega}\operatorname{curl}(\nabla q_E)\cdot\Phi-\nabla q_E\cdot\operatorname{curl}\Phi dx =$$ $$\int_{\Omega}q_{E}\operatorname{div}(\operatorname{curl}\Phi) dx-\int_{\partial\Omega}q_{E}\operatorname{curl}\Phi\cdot\nu ds 
			=0, 
		$$
		
		which implies that $\nabla q_E\times\nu=0\mathrm{~on~}\partial\Omega $, so we obtain,
		$$\Psi_E\times\nu=(\text{curl }\Phi_E)\times\nu=E\times\nu-\nabla q_E\times\nu=0\quad\text{on }\partial\Omega. $$ 
		
		Then $\Psi_E$ satisfy
		\begin{equation}
			\begin{cases}\operatorname{curl}\Psi_E=-\text{i}\omega\mu H+J_m-\operatorname{curl} G&\text{ in }\Omega,\\\mathrm{div}\Psi_E=0&\text{ in }\Omega. \\\Psi_E\times\nu=0&\text{ on }\partial\Omega,\end{cases}
		\end{equation}
		
		We have
		\begin{equation}\label{11}
			\|\Psi_E\|_{W^{1,p}(\Omega)}\leq C\|(H,J_m,\operatorname{curl}G)\|_{L^p(\Omega)},
		\end{equation}
		$$\operatorname{curl}\Phi_E= \Psi_E\in {W^{1,p}(\Omega)}. $$
		Set $\Psi_H=\operatorname{curl}\Phi_H$.  Observe that for every test function $\phi \in C^{\infty},$
		$$\int_\Omega\Psi_H\cdot\nabla\varphi+\varphi\mathrm{~div~}\Psi_Hdx=\int_\Omega\mathrm{curl~}\Phi_H\cdot\nabla\varphi dx$$
		$$=\int_\Omega\Phi_H\cdot\mathrm{curl~}\nabla\varphi dx+\int_{\partial\Omega}(\Psi_H \times   \nu) \cdot \nabla \phi
		ds=0,$$
		then   $\Psi_H \cdot \nu =0$ on $\partial\Omega$. 
		
		\begin{equation}
			\begin{cases}\operatorname{curl}\Psi_H=\text{i} \omega\varepsilon E+J_e&\text{ in }\Omega,\\\mathrm{div}\Psi_E=0&\text{ in }\Omega,\\\Psi_H\cdot\nu=0&\text{ on }\partial\Omega. \end{cases}
		\end{equation}
		
		We have
		\begin{equation}\label{13}
			\|\Psi_H\|_{W^{1,p}(\Omega)}\leq C\|(E,J_e)\|_{L^p(\Omega)},
		\end{equation}
		$$\operatorname{curl}\Phi_H= \Psi_E\in {W^{1,p}(\Omega)},$$

		since  $$\operatorname{curl}\Phi_H= \Psi_E\in {W^{1,p}(\Omega)},$$
		$$\operatorname{curl}\Phi_E= \Psi_E\in {W^{1,p}(\Omega)},$$
		$p \in(1,\frac{3}{2}) <3 $, then $\operatorname{curl}\Phi_H,\operatorname{curl}\Phi_E \in L^{p^*}({\Omega})$, and $p^{*}=\frac{3p}{3-p}\in(\frac{3}{2},3)<3$,
We will next use the other two sets of equations to estimate $\nabla q_E$ and $\nabla q_H$. 
		
		\begin{equation}\label{14}
			\begin{cases}-\mathrm{div}(\varepsilon \nabla q_E)=\mathrm{div}(\varepsilon \operatorname{curl}\Phi_E+\varepsilon G-i\omega^{-1}J_e)&\text{ in }\Omega,\\ \nabla q_E=0&\text{ on }\partial\Omega,\end{cases}
		\end{equation}
		\begin{equation}\label{15}
			\begin{cases}-\mathrm{div}(\mu \nabla q_H)=\mathrm{div}(\mu \operatorname{curl}\Phi_H+i\omega^{-1} \operatorname{curl}G+i\omega^{-1}J_m)&\text{ in }\Omega,\\ -(\mu \nabla q_H)\cdot \nu=(\mu \operatorname{curl}\Phi_H+i\omega^{-1} \operatorname{curl}G+i\omega^{-1}J_m)\cdot \nu&\text{ on }\partial\Omega,\end{cases}
		\end{equation}
		Therefore, by the $L^p$ theory for elliptics (Calderón–Zygmund estimate), we obtain $\nabla q_E,\nabla q_H \in L^{p^*}({\Omega})$,\\so $E=\nabla q_E+\operatorname{curl} \Phi_E \in L^{p^*}({\Omega})$ and $H=\nabla q_H+\operatorname{curl} \Phi_H \in L^{p^*}({\Omega})$
		and
		$$\|(E,H)\|_{L^{p^*}(\Omega)}\leq C(\|(E,H)\|_{L^p(\Omega)}+\|(J_e,J_m)\|_{L^{p^{*}}(\Omega)}+\|G\|_{W^{1,p^{*}}(\Omega)}). 
		$$
		
		\item $L^{p^{**}}$ estimate for $E,H$. 
		
		By \eqref{11} \eqref{13}, we obtain  
		$$\operatorname{curl} \Phi_E,\operatorname{curl} \Phi_H \in W^{1,p^{*}}({\Omega}) \hookrightarrow L^{p^{**}}({\Omega}),$$
		and 
		$$\|(\operatorname{curl} \Phi_E,\operatorname{curl} \Phi_H)\|_{W^{1,p^{*}}(\Omega)}\leq C\|(E,H)\|_{L^p{^{*}}(\Omega)}. $$
		\item
		by the $L^p$ throry for \eqref{14} \eqref{15} (Calderón–Zygmund estimate), we obtain $\nabla q_E,\nabla q_H \in L^{p^{**}}({\Omega})$,
		\\so $E=\nabla q_E+\operatorname{curl} \Phi_E \in L^{p^{**}}({\Omega})$ and $H=\nabla q_H+\operatorname{curl} \Phi_H \in L^{p^{**}}({\Omega})$,
		since sobolev embdding $p^{**}>3$, $\alpha=1-\frac{3}{p^{**}} \in(0,1)$,\\
		
		$\|(\operatorname{curl} \Phi_E,\operatorname{curl} \Phi_H)\|_ {C^{0,\alpha}(\bar{\Omega})}\leq C\|(\operatorname{curl} \Phi_E,\operatorname{curl} \Phi_H)\|_{W^{1,p^{**}}(\Omega)}\\ \leq C\|(E,H,J_e,J_m,\operatorname{curl}G)\|_{L^{p^{**}}(\Omega)}$. 
		
		\item$ C^{0,\alpha}$ estimate.  \\
		By applying Schauder esitmate for \eqref{14} \eqref{15}
		$$\nabla q_E,\nabla q_H \in C^{0,\alpha}(\overline{\Omega}),$$
		and
		\[
		\ \|\nabla q_E,\nabla q_H \|_{C^{0,\alpha}(\overline{\Omega})} \leq C(\|(\operatorname{curl} \Phi_E,\operatorname{curl} \Phi_H) \|_{C^{0,\alpha}(\overline{\Omega})} \] 
		\[+ \|(J_m,J_e) \|_{C^{0,\alpha}(\overline{\Omega})}+\|G \|_{C^{1,\alpha}(curl,\Omega)} ). \qedhere
		\]

	\end{enumerate}	
	
\end{proof}
\subsection{Proof of Proposition \ref{f}. }
\begin{proof}
	\begin{enumerate}[Step 1:]	
\item  Helmholtz decomposition.  There exist $q_E\in H_0^1(\Omega),q_H\in H^1(\Omega)\mathrm{~and~}\Phi_E,\Phi_H\in H^1(\Omega)\mathrm{~such~that}\\E=\nabla q_E+\operatorname{curl}\Phi_E\quad\operatorname{in} \Omega,\quad H=\nabla q_H+\operatorname{curl}\Phi_H\quad\operatorname{in} \Omega,$\\

$\begin{cases}\operatorname{div}\Phi_E=0\text{in }\Omega,\\\Phi_E\cdot\nu=0\text{ on }\partial\Omega,&\end{cases}\quad\begin{cases}\operatorname{div}\Phi_H=0\text{ in }\Omega,\\\Phi_H\times\nu=0\text{ on }\partial\Omega. &\end{cases}$

$$\|\Phi_E\|_{H^1(\Omega)}\leq C\left\|E\right\|_{L^2(\Omega)},\quad\|\Phi_H\|_{H^1(\Omega)}\leq C\left\|H\right\|_{L^2(\Omega)}. $$
\item$L^{p^{*}}$  for E. 

since 3. 1 step2 step3
	
$$\|\operatorname{curl}\Phi_E\|_{W^{1,p}(\Omega)}\leq C\|H\|_{L^p(\Omega)},$$
$$\|\operatorname{curl}\Phi_H\|_{W^{1,p}(\Omega)}\leq C\|E\|_{L^p(\Omega)},$$
$$E,H \in L^2(\Omega) \Rightarrow E,H \in L^p(\Omega),p\in (1,\frac{3}{2}),$$

$$\operatorname{curl} \Phi_E,\operatorname{curl} \Phi_H \in W^{1,p}(\Omega) \hookrightarrow L^{p^{*}}(\Omega),$$
$$p^{*}=\frac{3p}{3-p}\in(\frac{3}{2},3)<3. $$

	\begin{equation}
	\begin{cases}-\mathrm{div}(\varepsilon \nabla q_E)=\mathrm{div}(\varepsilon \operatorname{curl}\Phi_E)&\text{ in }\Omega,\\ \nabla q_E=0&\text{ on }\partial\Omega,\end{cases}
\end{equation}

we obtain $\nabla q_E(\Omega) \in L^{p^{*}}(\Omega) $
so $E \in L^{p^{*}}. $

$$\|\operatorname{curl}\Phi_H\|_{W^{1,p^{*}}(\Omega)}\leq C\|E\|_{L^{p^{*}}(\Omega)}. $$

\item Campanato estimate. . 
$$\operatorname{curl}\Phi_H \in W^{1,p^{*}}(\Omega) \hookrightarrow  L^{p^{**}}(\Omega) \hookrightarrow  L^{2,3\frac{p^{**}-2}{p^{**}}}(\Omega),$$
\begin{equation}
	\begin{cases}-\mathrm{div}(\mu \nabla q_H)=\mathrm{div}(\mu \operatorname{curl}\Phi_H)&\text{ in }\Omega,\\ -(\mu \nabla q_H)\cdot \nu=(\mu \operatorname{curl}\Phi_H)\cdot \nu&\text{ on }\partial\Omega,\end{cases}
\end{equation}

$$\nabla  q_H \in L^{2,min{(\gamma_0,\beta)}}(\Omega), $$ 

set $3\frac{p^{**}-2}{p^{**}}=\beta$ so $H \in  L^{2,min{(\gamma_0,\beta)}}(\Omega),$

since Lemma 5
$$\|\operatorname{curl}\Phi_E\|_{L^{2,min{(\gamma_0,\beta)+2}}(\Omega)}\leq C\|H\|_{L^{2,min{(\gamma_0,\beta)}}(\Omega)},$$
$$\operatorname{curl}\Phi_E \in L^{2,min{(\gamma_0,\beta)+2}}(\Omega)\cong C^{0,min{(\delta_0,1-\frac{3}{p^{**}})}}(\overline\Omega),  $$
Because $p$ is randomly selected in $(1,\frac{3}{2})$, So $p^{**}$  can arbitrarily approach infinity along with $p$.  Then  $min{(\delta_0,1-\frac{3}{p^{**}})}=\delta_0,$
$$\operatorname{curl}\Phi_E \in L^{2,min{(\gamma_0,\beta)+2}}(\Omega)\cong C^{0,\delta_0 }(\overline\Omega) . $$
\item  	By applying Schauder esitmate for  $\nabla q_E$
$$E \in C^{0,\alpha}(\overline\Omega), \quad \alpha \in (0,\delta_{0}] . $$  
\end{enumerate}
\end{proof}
 
\subsection{Proof of Theorem \ref{d}. }
We only give the differences between non-homogeneous and homogeneous. 
\begin{proof}
	\begin{enumerate}[Step 1:]	
		\item  Helmholtz decomposition. 
		$$E-G=\nabla q_E+\operatorname{curl}\Phi_E\quad\operatorname{in} \Omega,\quad H=\nabla q_H+\operatorname{curl}\Phi_H\quad\operatorname{in} \Omega,$$
	\item $L^{p^{*}}$ estimate.  \\
	  Estimate for $\operatorname{curl}\Phi_E \operatorname{curl}\Phi_H,$
		$$\|\operatorname{curl}\Phi_E\|_{W^{1,p}(\Omega)}\leq C\|(H,J_m,\operatorname{curl}G)\|_{L^p(\Omega)},$$
	   $$\|\operatorname{curl}\Phi_H\|_{W^{1,p}(\Omega)}\leq C\|(E,J_e\|_{L^p(\Omega)}. $$
	 Estimate ($\nabla q_E$). 
	\begin{equation*}
		\begin{cases}
		-\mathrm{div}(\varepsilon \nabla q_E)=\mathrm{div}(\varepsilon \operatorname{curl}\Phi_E+\varepsilon G-i\omega^{-1}J_e)&\text{ in }\Omega,\\ \nabla q_E=0&\text{ on }\partial\Omega,\end{cases}
	\end{equation*}
	
	$$\|\nabla q_E\|_{L^{p^{*}}(\Omega)}\leq C\|(\operatorname{curl}\Phi_E,J_e,G)\|_{L^{p^{*}}(\Omega)},$$
	
$$E \in L^{p^{*}}(\Omega). $$
	
	 $W^{1,p^{*}}$ Estimate. 
	
	$$\|\operatorname{curl}\Phi_H\|_{W^{1,p^{*}}(\Omega)}\leq C\|E,J_e\|_{L^{p^{*}}(\Omega)}. $$

	\item Campanato estimate . 
	$$\operatorname{curl}\Phi_H \in W^{1,p^{*}}(\Omega) \hookrightarrow  L^{p^{**}}(\Omega) \hookrightarrow  L^{2,3\frac{p^{**}-2}{p^{**}}} (\Omega), $$
	\begin{equation*}
		\begin{cases}-\mathrm{div}(\mu \nabla q_H)=\mathrm{div}(\mu \operatorname{curl}\Phi_H+i\omega^{-1} \operatorname{curl}G+i\omega^{-1}J_m)&\text{ in }\Omega,\\ -(\mu \nabla q_H)\cdot \nu=(\mu \operatorname{curl}\Phi_H+i\omega^{-1} \operatorname{curl}G+i\omega^{-1}J_m)\cdot \nu&\text{ on }\partial\Omega,\end{cases}
	\end{equation*}
	need  $\operatorname{curl}G,J_m \in L^{2,\gamma_0}$ 
	can obtain $\nabla  q_H \in L^{2,min{(\gamma_0,\beta)}}(\Omega) $ set $3\frac{p^{**}-2}{p^{**}}=\beta$ so $H \in  L^{2,min{(\gamma_0,\beta)}}(\Omega),$
	$$\|\operatorname{curl}\Phi_E\|_{L^{2,min{(\gamma_0,\beta)+2}}(\Omega)}\leq C\|H,J_m,\operatorname{curl}G\|_{L^{2,min{(\gamma_0,\beta)}}(\Omega)},$$
	$$\operatorname{curl}\Phi_E \in L^{2,min{(\gamma_0,\beta)+2}}(\Omega)\cong C^{0,min{(\delta_0,1-\frac{3}{p^{**}})}}(\overline{\Omega})  . $$
	Because $p$ is randomly selected in $(1,\frac{3}{2})$, So $p^{**}$  can arbitrarily approach infinity along with $p$.  Then  $min{(\delta_0,1-\frac{3}{p^{**}})}=\delta_0,$
	$$\operatorname{curl}\Phi_E \in L^{2,min{(\gamma_0,\beta)+2}}(\Omega)\cong C^{0,\delta_0 }(\overline{\Omega}) . $$
	\item  	By applying Schauder esitmate for  $\nabla q_E,$
	$$E \in C^{0,\alpha}(\overline{\Omega}),\quad \alpha \in (0,\delta_{0}]. $$

    From the proof above, we can see that we need $G,J_e \in C^{0,\alpha}$,
	$\operatorname{curl}G,J_m\in   L^{2,\gamma_0}. $
	\end{enumerate}

\end{proof}
\subsection{Proof of Theorem \ref{b}. }
Based on the comparison between homogeneous and non-homogeneous above, we will now only prove the non-homogeneous case. 
In the  Theorem \ref{b}, $\varepsilon(x),\mu(x) \in C^{0} (\bar{\Omega},\mathbb{R}^3)$, There are many steps that are the same as Theorem \ref{m}, so we will directly use it. 

\begin{proof}
\begin{enumerate}[Step 1:]	
	\item  Helmholtz decomposition. There exist $q_E\in H_0^1(\Omega),q_H\in H^1(\Omega)\\
	\mathrm{~and~}\Phi_E,\Phi_H\in H^1(\Omega)\mathrm{~such~that}\\E-G=\nabla q_E+\operatorname{curl}\Phi_E\quad\operatorname{in} \Omega,\quad H=\nabla q_H+\operatorname{curl}\Phi_H\quad\operatorname{in},$\\

	$\begin{cases}\operatorname{div}\Phi_E=0\text{ in }\Omega,\\\Phi_E\cdot\nu=0\text{ on }\partial\Omega,&\end{cases}\quad\begin{cases}\operatorname{div}\Phi_H=0\text{ in }\Omega,\\\Phi_H\times\nu=0\text{ on }\partial\Omega. &\end{cases}$
	
	$$\|\Phi_E\|_{H^1(\Omega)}\leq C\left\|E\right\|_{L^2(\Omega)},\quad\|\Phi_H\|_{H^1(\Omega)}\leq C\left\|H\right\|_{L^2(\Omega)}. $$
	
	\item $L^{p^{*}}$ estimate for E,H .

	Set $\Psi_E=\operatorname{curl}\Phi_E$.  Observe that for every test function $\phi \in C^{\infty},$
	
	since integration by parts  \\
	$$
		\int_{\Omega}\operatorname{curl}(\nabla q_E)\cdot\Phi-\nabla q_E\cdot\operatorname{curl}\Phi dx$$ $$=\int_{\Omega}q_{E}\operatorname{div}(\operatorname{curl}\Phi) dx-\int_{\partial\Omega}q_{E}\operatorname{curl}\Phi\cdot\nu ds \\
		=0, 
	$$
	
	which implies that $\nabla q_E\times\nu=0\mathrm{~on~}\partial\Omega $, so we obtain 
	$$\Psi_E\times\nu=(\text{curl }\Phi_E)\times\nu=E\times\nu-\nabla q_E\times\nu=0\quad\text{on }\partial\Omega,$$ 
	
	Then $\Psi_E$ satisfy
	\begin{equation}
		\begin{cases}\operatorname{curl}\Psi_E=-\text{i}\omega\mu H+J_m-\operatorname{curl} G&\text{ in }\Omega,\\\mathrm{div}\Psi_E=0&\text{ in }\Omega,\\\Psi_E\times\nu=0&\text{ on }\partial\Omega. \end{cases}
	\end{equation}
	
	We have
	\begin{equation}\label{19}
		\|\Psi_E\|_{W^{1,p}(\Omega)}\leq C\|(H,J_m,\operatorname{curl}G)\|_{L^p(\Omega)},
	\end{equation}
	$$\operatorname{curl}\Phi_E= \Psi_E\in {W^{1,p}(\Omega)},$$
	Set $\Psi_H=\operatorname{curl}\Phi_H$. Observe that for every test function $\phi \in C^{\infty},$
	$$\int_\Omega\Psi_H\cdot\nabla\varphi+\varphi\mathrm{~div~}\Psi_Hdx=\int_\Omega\mathrm{curl~}\Phi_H\cdot\nabla\varphi dx$$
	$$=\int_\Omega\Phi_H\cdot\mathrm{curl~}\nabla\varphi dx+\int_{\partial\Omega}(\Psi_H \times   \nu) \cdot \nabla \phi
	ds=0,$$
	then   $\Psi_H \cdot \nu =0$ on $\partial\Omega. $
	
	\begin{equation}
		\begin{cases}\operatorname{curl}\Psi_H=\text{i} \omega\varepsilon E+J_e&\text{ in }\Omega,\\\mathrm{div}\Psi_E=0&\text{ in }\Omega,\\\Psi_H\cdot\nu=0&\text{ on }\partial\Omega. \end{cases}
	\end{equation}
	
	We have
	\begin{equation}\label{21}
		\|\Psi_H\|_{W^{1,p}(\Omega)}\leq C\|(E,J_e)\|_{L^p(\Omega)},
	\end{equation}
	$$\operatorname{curl}\Phi_H= \Psi_E\in {W^{1,p}(\Omega)},$$
	
	since  $$\operatorname{curl}\Phi_H= \Psi_E\in {W^{1,p}(\Omega)},$$
	$$\operatorname{curl}\Phi_E= \Psi_E\in {W^{1,p}(\Omega)},$$
	$p \in(1,\frac{3}{2}) <3 $, then $\operatorname{curl}\Phi_H,\operatorname{curl}\Phi_E \in L^{p^*}$,and $p^{*}=\frac{3p}{3-p}\in(\frac{3}{2},3)<3. $
	We will next use the other two sets of equations to estimate $\nabla q_E$ and $\nabla q_H$. 
	
	\begin{equation}\label{22}
		\begin{cases}-\mathrm{div}(\varepsilon \nabla q_E)=\mathrm{div}(\varepsilon \operatorname{curl}\Phi_E+\varepsilon G-i\omega^{-1}J_e)&\text{ in }\Omega,\\ \nabla q_E=0&\text{ on }\partial\Omega,\end{cases}
	\end{equation}
	\begin{equation}\label{23}
		\begin{cases}-\mathrm{div}(\mu \nabla q_H)=\mathrm{div}(\mu \operatorname{curl}\Phi_H+i\omega^{-1} \operatorname{curl}G+i\omega^{-1}J_m)&\text{ in }\Omega,\\ -(\mu \nabla q_H)\cdot \nu=(\mu \operatorname{curl}\Phi_H+i\omega^{-1} \operatorname{curl}G+i\omega^{-1}J_m)\cdot \nu&\text{ on }\partial\Omega,\end{cases}
	\end{equation}
	Therefore,by the $L^p$ theory for elliptics (Calderón–Zygmund estimate), we obtain $\nabla q_E,\nabla q_H \in L^{p^*}(\Omega)$,\\so $E=\nabla q_E+\operatorname{curl} \Phi_E \in L^{p^*}(\Omega)$ and $H=\nabla q_H+\operatorname{curl} \Phi_H \in L^{p^*}(\Omega)$
	and
	$$\|(E,H)\|_{L^{p^*}(\Omega)}\leq C(\|(E,H)\|_{L^p(\Omega)}+\|(J_e,J_m)\|_{L^{p^{*}}(\Omega)}+\|G\|_{W^{1,p^{*}}(\Omega)}). 
	$$
	
	\item $L^{p^{**}}$ estimate for $E,H$. 
	
	By \eqref{19} \eqref{21}, we obtain  
	$$\operatorname{curl} \Phi_E,\operatorname{curl} \Phi_H \in W^{1,p^{*}}(\Omega) \hookrightarrow L^{p^{**}}(\Omega),$$
	and 
	$$\|(\operatorname{curl} \Phi_E,\operatorname{curl} \Phi_H)\|_{W^{1,p^{*}}(\Omega)}\leq C\|(E,H,J_e,J_m,\operatorname{curl}  G)\|_{L^p{^{*}}(\Omega)}. $$
	\item
	By the $L^p$ theory for \eqref{22} \eqref{23} (Calderón–Zygmund estimate), we obtain $\nabla q_E,\nabla q_H \in L^{min (q,{p^{**}})}=L^q({\Omega}),$
	\\so $E=\nabla q_E+\operatorname{curl} \Phi_E \in L^q(\overline{\Omega})$ and $H=\nabla q_H+\operatorname{curl} \Phi_H \in L^q(\overline{\Omega}). $ \\
	Finally $$ L^q({\Omega}) \hookrightarrow \mathcal{L}^{2,3\frac{q-2}{q}}({\Omega}). $$
\end{enumerate}	
\end{proof}
\subsection{Proof of Theorem \ref{c}. }

The proof process is exactly the same as the first four steps of Theorem \ref{m}.  In the  step 5, since the Schauder estimate.  \( \varepsilon  \) is \( C^{0,\alpha} \)  but \( \mu \) is \( C^0 \), we can only conclude that \( E \) is \( C^{0,\alpha} \).

\appendix
\section{ Appendix }
In the appendix, we present several key functional spaces and their properties utilized in this paper.  Since most of these spaces and properties are widely recognized. \cite{CHEN}

\begin{definition}[Morrey Space]
	We denote by $L^{p,\lambda}( \Omega )$ the Morrey spaces, and for $1\leq p<\infty$ and
	$\lambda\in(0,n+p)$, its norm is defined as
	
	$$\|u\|_{{L}^{p,\lambda}}(\Omega)=\{\sup_{x\in\Omega,0<r<d}r^{-\lambda}\int_{\Omega(x,r)}|u(z)|^pdz\}^{\frac1p},$$
	
	where $\Omega(x,r)=\Omega\cap B(x,r),B(x,r)$ is a ball of which centre is $x\in \Omega . $ We call $u(x) \in  L^{p,\lambda}(\Omega)$ if $\|u\|_{L^{p,\lambda}}(\Omega)<\infty. $
\end{definition}

\begin{definition}[Campanato Space]
We denote by $L ^{p, \lambda }( \Omega )$ the Campanato spaces, and for the same $p,\lambda$
used in Definition 11, the norm of the Campanato spaces is expressed as
$$\|u\|_{\mathcal{L}^{p,\lambda}}(\Omega)=\|u\|_{{L}^{p,\lambda}}(\Omega)+\{\sup_{x\in\Omega,0<r<d}r^{-\lambda}\int_{\Omega(x,r)}|u(z)-\bar{u}_{x,r}|^pdz\}^{\frac{1}{p}},$$

where $\bar{u}_{x,r}=\frac1{|\Omega(x,r)|}\int_{\Omega(x,r)}u(z)dz. $

The Morrey space and Campanato space have the following two well-known and important properties. 
\end{definition}
In the properties, we always assume $ n = 3 $. 
\begin{proposition}
	If  $\lambda \in (3,5) $, then  $ \mathcal{L}^{2,\lambda}({\Omega}) \cong C^{0,\frac{\lambda-3}2}(\overline{\Omega})$.  
\end{proposition}
\begin{proposition}
	$\mathcal{L}^{p,\lambda}(\Omega)\cong{L}^{p,\lambda}$,  for $0\leq\lambda<n$. 
\end{proposition}
\begin{proposition}
	The embedding $L^p(\Omega)\hookrightarrow \mathcal{L}^{2,3\frac{p-2}p}(\Omega)$ is continuous. 
\end{proposition}
\begin{proposition}
	If $\lambda<3,u\in L^{2}(\Omega)$ and $\nabla u\in \mathcal{L}^{2,\lambda}(\Omega)$ then $u\in \mathcal{L}^{2,2+\lambda}(\Omega)$  and the embedding is continuous. 
\end{proposition}

\end{document}